\documentstyle[12pt]{article}
\input epsf
\def\del#1{}
\def\dell#1{}
\def\ins#1{#1}

\newcommand{\limfunc}[1]{{\rm \mathop{#1}}}
\newcommand{\text}[1]{{\mbox {#1}}}
\newcommand{\func}[1]{{\rm #1} \,}
\newtheorem{theorem}{Theorem}
\newtheorem{lemma}{Lemma}
\newtheorem{proposition}{Proposition}
\newtheorem{remark}{Remark}
\newtheorem{corollary}{Corollary}
\newtheorem{definition}{Definition}

\newlength{\rig}
\newlength{\rigg}
\newlength{\hei}
\newcommand{\fgr}[3]{
\setlength{\hei}{7.8cm}
\setlength{\rig}{0.1\textwidth}
\setlength{\rigg}{0.8\textwidth}
\begin{figure}
\rule{\rig}{0in}
\epsfxsize=13cm
\epsffile[1 1 953 638]{#1}
\caption{#3}\label{#2}
\end{figure}
}

\begin{document}
\setcounter{page}{1}

\author{
{\bf Anton Savin\footnotemark\rule{5pt}{0pt}}\\
Moscow State University\\
e-mail: antonsavin@mtu-net.ru\\
\&\\
{\addtocounter{footnote}{-1}
{\bf Boris Sternin}\footnotemark}%
\addtocounter{footnote}{-1}\thanks{
Supported by Institut f\"ur Mathematik, Universit\"at Potsdam,
and RFBR grants
Nos.~97-01-00703, 97-02-16722a, 99-01-01254, 99-01-01100, and
Soros Foundation grant a99-1763.}\\
Moscow State University\\
e-mail: sternine@mtu-net.ru\\
\&\\
{\bf Bert-Wolfgang Schulze}\\
Potsdam University\\
e-mail: schulze@math.uni-potsdam.de\\
}
\title{\bf
\ins{The} Homotopy Classification and \ins{the}
Index of Boundary Value Problems
for General Elliptic Operators
}
\date{October, 1999}

\maketitle

\begin{abstract}

We give the homotopy classification and compute \ins{the} index of boundary
value problems for elliptic equations. The classical case of operators
\del{which}\ins{that}  satisfy \ins{the} Atiyah--Bott condition
is studied first. We \ins{also} consider\del{ also}
the general case of boundary value problems for operators that
\del{possibly
violate}\ins{do not necessarily satisfy the} Atiyah--Bott condition.
\end{abstract}

{\bf Keywords}:
elliptic boundary value problems, Atiyah--Bott condition, index theory,
$K$-theory, homotopy classification

\vspace{15mm}

{\bf 1991 AMS classification}: Primary 58G20, Secondary 58G10, 19M05, 35J55

\vfill

\newpage

\tableofcontents

\section*{Introduction}
\addcontentsline{toc}{section}{Introduction}

The present paper \del{is devoted to}\ins{deals with}
the theory of boundary value
problems for elliptic equations. The well-known Shapiro--Lopatinskii
condition (\del{see }e.g.\ins{, see}~\cite{Hor3})
describes the class of elliptic boundary
value problems, i.e.\ins{,} problems defining Fredholm operators in Sobolev
spaces. On the other hand, this condition represents the obstruction
to \ins{the} existence of \del{correct}\ins{well-posed}
(Fredholm) boundary value problems for
elliptic operators on manifolds with boundary. Moreover, from the
topological point of view, \ins{the}
Shapiro--Lopatinskii condition guarantees \ins{the existence of}
a homotopy of the homogeneous principal symbol of the given
elliptic
operator to a
symbol independent of the cotangent variables in
\del{the}\ins{a}  neighborhood of the
boundary. This re\del{formulation}\ins{statement}
of the Shapiro--Lopatinskii condition
was found by Atiyah and Bott\del{ in} \cite{AtBo2}\dell{,}
\del{therefore this condition}\ins{and}
is \del{frequently}\ins{often}  called \ins{the}
Atiyah--Bott condition.  From the analytical
point of view, \ins{the} Shapiro--Lopatinskii condition
\del{makes it possible}\ins{permits one}  to reduce
a boundary value problem to \del{an}\ins{a zero-order}
elliptic operator\del{ of order zero, such} that \ins{is
a bundle homomorphism} in
a neighborhood of the boundary\ins{.}\del{
it is an operator of multiplication.}
A reduction of this kind is fundamental for the homotopy
classification of elliptic boundary value problems\dell{,}
and for the derivation
of the corresponding index formula.

As\del{ it} was \del{emphasized}\ins{already pointed out},
not all operators on a manifold with
boundary admit  \del{correct}\ins{well-posed classical}
boundary value problems.
The theory of boundary value problems for general elliptic operators
(which need not satisfy \ins{the} Atiyah--Bott condition) was
constructed in \cite{ScSS18}\ins{.}\footnote{In this
connection\ins{,} we\del{ would like to}
note \del{an}\ins{the} book\del{ of
Booss--Wojciechowski} \cite{BBW1} \ins{by Booss and
Wojciechowski}, where, in
particular, the theory of boundary value problems is constructed for
operators similar to the Dirac operator, which is an important
geometric\del{al} operator
that does not satisfy the Atiyah--Bott condition.}\dell{.}
In this theory\ins{,} \ins{the}  violation of the
Shapiro--Lopatinskii--Atiyah--Bott condition does not allow
one to reduce an elliptic boundary value problem to
a zero-order operator. Nonetheless,
\del{such}\ins{the} reduction \del{becomes}\ins{is}
possible\dell{,}
\del{when}\ins{if}
the boundary value problem possesses \ins{certain}
symmetries. In this case\ins{,} it
is also possible to give \del{the}\ins{a}
homotopy classification of elliptic boundary
value problems with symmetries and\del{ to} obtain
\del{the}\ins{an}  index
formula.\vspace{0.1cm}

Let us describe the contents of the paper in more detail.

\vspace{0.1cm}

The first part of the paper \del{is devoted to}\ins{deals with}
classical boundary value problems.
The main result here is the homotopy classification of boundary
value problems. More precisely, it is shown that a boundary value problem
satisfying the Shapiro--Lopatinskii condition admits a reduction to a
zero-order operator that requires no boundary conditions,
i.e.\ins{,} there is an
isomorphism
\[
\varepsilon _m:\limfunc{Ell}^m\left( M\right) \longrightarrow
\limfunc{Ell}^0\left( M\right)
\]
of the group of stably homotopic elliptic boundary value problems
for operators of order $m\geq 1$ and \del{an analogous}\ins{a
similar}  group for zero-order
operators. Furthermore, the operators of order zero are
classified, \ins{in the same way} as
in elliptic theory on closed manifolds, by their principal symbols
\[
\chi :\limfunc{Ell}^0\left( M\right) \longrightarrow K\left(%
T^{*}\left( M\backslash \partial M\right) \right) ,
\]
where $K\left( T^{*}\left( M\backslash
\partial M\right) \right) $ \del{denotes}\ins{is the}
$K$-group with compact support\ins{s}.
In this section of \ins{the paper }%
we follow H\"ormander \cite{Hor3}, \del{where}\ins{who realized}  the
topological method \del{of}\ins{due to}
Atiyah--Bott \cite{AtBo2} \del{is realized in}\ins{by}  explicit
homotopies of boundary value problems.\del{Another }
\del{characteristic feature of
this approach is}\ins{We point out} that
\del{the}\ins{H\"ormander's}  homotopies of classical boundary value
problems \del{are carried out without the}\ins{do not} use
\del{of}\ins{the}  complete Boutet de Monvel algebra
\cite{Bout2}. This \del{makes it possible}\ins{permits one}
to obtain\dell{:} \ins{the} homotopy
classification of boundary value problems\dell{,} \ins{and}
the corresponding index formula\dell{,}
and \del{at the same time to}\ins{simultaneously}
prove \ins{the} Atiyah--Bott theorem on the obstruction
to \ins{the} existence of classical boundary value problems.

\del{T}\ins{In t}he second part of the paper\ins{,}
\del{is devoted to}\ins{we consider}  boundary value problems
\cite{ScSS18} for operators \del{which violate}\ins{that do
not satisfy the}
Shapiro--Lopatinskii--Atiyah--Bott condition.
The\ins{se} boundary value problems \del{under consideration
are of}\ins{have} the form
\begin{equation}
\left\{
\begin{array}{ll}
Du=f , & u\in H^s\left( M,E\right) ,\quad f\in H^{s-m}\left( M,F\right) ,
\vspace{0.2cm}\\
Bj_{\partial M}^{m-1}u=g,\quad & g\in \func{Im}P\subset H^\sigma
\left( \partial M,G\right) ,
\end{array}
\right.  \label{one}
\end{equation}
where $B$ is a boundary operator with\del{ its} range
\ins{contained} in \del{a subspace}\ins{the range} $\func{Im}P$
of \ins{a pseudodifferential projection $P$ in}
the Sobolev space on the boundary of $M$\del{, which is}
\del{determined by a pseudodifferential projection $P$, }\ins{and}
$j_{\partial M}^{m-1}$
is \del{a}\ins{the}
composition of the jet of order $m-1$ and the restriction to the
boundary $\partial M$. \del{Because}\ins{In view}  of the Atiyah--Bott
obstruction\ins{,} the boundary value problem
(\ref{odin}) can\del{ }not
be reduced to an operator of order zero.
However, an arbitrary elliptic boundary value problem
\del{is}\ins{can be}  reduced in this
case to the so-called {\em spectral boundary value problem}
\cite{APS1,NScSS3} for a first-order operator.

Further simplification of the boundary value problem is possible under
additional assumptions on the subspace defined by
\del{a}\ins{the}  pseudodifferential
projection \del{from}\ins{on}  the right-hand side
\del{of}\ins{in} (\ref{one}).
An example of such assumptions is given by\del{ the} parity conditions imposed
on the principal symbol of \ins{the} projection
(see \cite{SaSt1,SaSt2}). Precise definitions will be given below,
\del{but}\ins{and}
for now \del{let us}\ins{we only}
mention that these conditions can be reformulated as
conditions under which the operator $D$ of the boundary value problem
extends to the double of the manifold with boundary in a symmetric way.
\del{Such reformulation}\ins{This restatement}
shows that the parity condition is a generalization
of \ins{the} Atiyah--Bott condition, which guarantees
\ins{the existence of} a homotopy of the principal
symbol of the operator to \ins{the}
identity symbol in a neighborhood of the
boundary (and, of course, the possibility of extension to the double).

Under the above-mentioned parity conditions\ins{,} \ins{the}
stable homotopy classification
modulo 2-torsion is obtained for elliptic boundary value
problems\ins{.} \ins{It has the form}
\begin{equation}
\limfunc{Ell}^{ev/odd}\left( M,\partial M\right) \otimes {\bf Z}%
\left[ \frac 12\right] \simeq K\left( T^{*}\left( M\backslash \partial
M\right) \right) \otimes {\bf Z}\left[ \frac 12\right] \oplus {\bf Z}%
\left[ \frac 12\right] ,  \label{xstar}
\end{equation}
where the\del{ map on the} first term in the sum is
\del{defined}\ins{determined}  by the principal
symbol of \ins{the} boundary value problem\dell{,}
\del{while}\ins{and}  the second component is
given by the value of \del{the}\ins{a} functional $d$ \ins{on
the set} of subspaces\dell{,} defined by
pseudodifferential projections \ins{(this functional was
defined in}
\cite{SaSt1,SaSt2}\ins{)}.
The functional $d$ \ins{is} equal\del{s} \ins{to}
the \ins{Atiyah--Patodi--Singer} spectral $\eta $-invariant
\del{of Atiyah--Patodi--Singer }\cite{APS1} of a\ins{n
admissible} self-adjoint operator \del{in}\ins{for} the
case\dell{,}
\del{when}\ins{in which} the subspace \ins{in question}
is\del{ defined as} the nonnegative spectral subspace \del{for an
admissible self-adjoint}\ins{of that}  operator \cite{Gil7,SaSt1}.

The index formula for boundary value problems with parity conditions
\cite{SaSt1,SaSt2}\del{ is in a} \del{natural way an immediate
consequence of}\ins{readily follows from}
the homotopy classification (\ref{xstar}) \ins{in a natural way}.

The authors are grateful to V.~E.~Nazaikinskii for numerous
useful discussions. The results of the paper were reported at
\ins{the}
international conference "Workshop in Partial Differential
Equations\ins{,}"\dell{,}
July 1999, Potsdam, Germany\dell{,} and also at the conference "Jean
Leray 1999\ins{,}"\dell{,}
August 1999, Karlskrona-Ronneby, Sweden.

\section{Classical boundary value problems}

\subsection{Basic definitions}

Let $M$ be a compact smooth manifold with boundary $X=\partial
M$ and%
\del{ $D$ be
an elliptic} \del{differential operator on $M$ of order $m\geq 1$}
\[
D:H^s\left( M,E\right) \longrightarrow H^{s-m}\left( M,F\right)
\]
\ins{an elliptic differential operator of order $m\geq 1$}
acting in Sobolev spaces of sections of vector bundles $E,F$
\ins{over $M$.}
\ins{The} \del{O}\ins{o}perator $D$ is not\del{ a}
Fredholm\ins{,}\del{ operator in Sobolev} \del{spaces:}\ins{since}  its
kernel $\ker D$ is\del{ an} infinite-dimensional\del{ space}.
To define a Fredholm
operator, let us \del{add to}\ins{equip}  $D$ \ins{with}
some boundary conditions.\del{ They are defined in}
\del{the following way. }\ins{To this end,} \del{W}\ins{w}e
\del{fix}\ins{choose}  a collar neighborhood
$U=X\times \left[0,1\right) \subset M$ of\del{ the boundary} $X$
with normal coordinate $t\in \left[0,1\right)$. Consider jets
of order $m-1$ in the normal direction to the boundary composed
with the operator of restriction to the boundary
%\[
%j_X^{m-1}u=\left( \left. u\right| _X,\left. -i\frac \partial {\partial
%t}u\right| _X,\ldots ,\left. \left( -i\frac \partial {\partial t}\right)
%^{m-1}u\right| _X\right) :H^s\left( M,E\right) \rightarrow
%\bigoplus_{k=0}^{m-1}H^{s-1/2-k}\left( X,\left. E\right| _X\right) .
%\]
\[
\begin{array}{c}
j_X^{m-1}u=\left( \left. u\right| _X,\left. -i\frac \partial {\partial
t}u\right| _X,\ldots ,\left. \left( -i\frac \partial {\partial t}\right)
^{m-1}u\right| _X\right) , \\
j_X^{m-1}:H^s\left( M,E\right) \rightarrow
\bigoplus_{k=0}^{m-1}H^{s-1/2-k}\left( X,\left. E\right| _X\right) .
\end{array}
\]
A {\em classical boundary value problem} for\del{ elliptic differential
operator} $D$ is a system of equations of the form
\begin{equation}
\left\{
\begin{array}{ll}
Du=f, & u\in H^s\left( M,E\right) ,f\in H^{s-m}\left( M,F\right) ,
\vspace{0.2cm}  \\
Bj_X^{m-1}u=g, & g\in H^\sigma \left( X,G\right) ,
\end{array}
\right.  \label{star1}
\end{equation}
where
\begin{equation}
B:\bigoplus_{k=0}^{m-1}H^{s-1/2-k}\left( X,\left. E\right| _X\right)
\longrightarrow H^\sigma \left( X,G\right)  \label{mu1}
\end{equation}
is a pseudodifferential operator on the boundary;
\del{in addition, }the orders of its components and the indices of Sobolev
spaces in (\ref{mu1}) are supposed to be compatible in a natural way
(\del{see }e.g.\ins{,} \ins{see}~\cite{Hor3}). For
brevity\ins{,} \ins{the}  boundary value problem
$\left( D,B\right) $ will\del{ be} sometimes \ins{be}
denoted by ${\cal D}.$

On the cotangent sphere bundle $S^{*}X$ of $X$\ins{,} \ins{we}
consider \del{a}\ins{the}  vector bundle
\[
L_{+}\left( D\right) \subset \pi ^{*}E^m,\qquad \pi
:S^{*}X\rightarrow X,
\]
\ins{whose fiber over a point $(x,\xi)\in S^*X$ is}
\del{defined as }the \del{Cauchy}\ins{subspace of}  initial data
\del{subspace for}\ins{of bounded} solutions of
\del{an}\ins{the}  ordinary
differential equation
\[
\sigma \left( D\right) \left( x,0,\xi ^{\prime },-i\frac d{dt}\right)
u\left( t\right) =0,\qquad \left( x,\xi ^{\prime }\right) \in S^{*}X,
\]
\ins{with constant coefficients
on the half-line $\left\{t\geq 0\right\}$.}\del{
which are bounded as  $t\rightarrow +\infty .$}
The complementary subbundle\dell{,}
corresponding to solutions bounded as $t\rightarrow
-\infty$\dell{,} is denoted by
$L_{-}\left( D\right)$. The subbundles $L_{\pm }\left(
D\right) $ are\dell{,}
obviously\dell{,} determined by the restriction of
the principal symbol of $D$ to the boundary.

The restriction
\begin{equation}
L_{+}\left( D\right) \stackrel{\sigma \left( B\right) }{\longrightarrow }\pi
^{*}G  \label{alpha1}
\end{equation}
of the principal symbol of \ins{the}
boundary operator $B$ to the
subbundle $L_{+}\left( D\right) \subset\pi^{*}E^m$ is called
\ins{the} {\em boundary
symbol} of classical boundary value problem $\left( D,B\right)
$\ins{.}

\ins{The} \del{B}\ins{b}oundary value problem
$\left(D,B\right) $ is \del{called}\ins{said to be}
{\em elliptic}\dell{,}
if its boundary symbol is an isomorphism of vector bundles.

\begin{proposition}
{\em (\del{see }e.g.\ins{,} \ins{see}~\cite{Hor3})} \ins{The}
\del{B}\ins{b}oundary value problem {\em (\ref{star1})}
has \ins{the} Fredholm property if and only if it is elliptic.
\end{proposition}

The ellipticity condition (\ref{alpha1}) imposes an essential
restriction on the bundle $L_{+}\left( D\right) $: for the existence
of an elliptic boundary value problem for $D$, \ins{it is
necessary that}  this bundle
\del{has to }be isomorphic to a bundle\dell{,} lifted
from\del{ the boundary}
$X$\dell{,}\ins{;} \del{while}\ins{the choice of}
a specific lifting (\ref{alpha1})
\del{defines}\ins{determines}  the boundary conditions.
Atiyah and Bott \cite{AtBo2} noted that this condition can be
re\del{formul}\ins{st}ated
in terms of the principal symbol of\del{ operator} $D$ in the following form:
the restriction of\del{ the symbol}
$\sigma \left(D\right) $ to\del{ the boundary}
\del{of the manifold }\ins{$X$}
is stably homotopic to the symbol of \ins{a} multiplication operator,
that is\ins{,}
\begin{equation}
\sigma \left( D\right) \left( x,0,\xi \right) \sim \sigma ^{\prime }\left(
x\right) ,  \label{sstar1}
\end{equation}
or\ins{,} in terms of $K$-theory\ins{,}
\[
\left[ \sigma \left( D\right) \right] \in \func{Im}\left\{ K\left(
T^{*}\left( M\backslash \partial M\right) \right) \stackrel{\pi ^{*}}{%
\longrightarrow }K\left( T^{*}M\right) \right\} ,\quad \pi
:T^{*}M\rightarrow \left. T^{*}M\right/ \left. T^{*}M\right| _X.
\]
\del{In addition}\ins{Furthermore}, the choice  of
a boundary condition \del{defines}\ins{determines} a certain homotopy
of the form (\ref{sstar1}), which
\del{defines}\ins{specifies}  an element
\[
\left[ \sigma \left( D,B\right) \right] \in K\left( T^{*}\left( M\backslash
\partial M\right) \right) .
\]
It turns out (see \del{s}\ins{S}ection 2) that this element classifies
the boundary value problem $\left(D,B\right) $ up to stable
homotopy \ins{equivalence}.

In the next section\ins{,}
we carry out the homotopy classification of elliptic
boundary value problems. To this end, \del{it is
necessary}\ins{we have}  to enlarge the class
of operators\dell{,} for which\del{ we can pose} boundary value problems
\ins{will be posed}. Namely,
we deal with elliptic operators $D$ on a manifold \ins{$M$}
with boundary \ins{$X$} which
satisfy the following conditions.

\begin{enumerate}
\item  \del{The operator $D$ i}\ins{I}n a small neighborhood
\del{of the boundary }$X\times\left[ 0,2\varepsilon \right) $
\ins{of the boundary,}  \del{is of}\ins{$D$ has} the form
\begin{equation}
D=\sum_{k=0}^mD_k\left( t\right) \left( -i\frac \partial {\partial t}\right)
^{m-k},  \label{omeg}
\end{equation}
where \ins{the}
$\!D_k\!\left( t\right)$ are smooth families of pseudodifferential
operators on\del{ the boundary} $\!X\!$, $\limfunc{ord}D_k\!\left( t\right)\! =k,$
such that\del{ the family} $D_0\left( t\right) $ consists of isomorphisms
of vector bundles. By identifying the vector bundles
in which
$D$ acts with the help of this isomorphism in a neighborhood
of the boundary, \del{it}\ins{we} can\del{ be} assume\del{d} that the
coefficient $D_0\left(t\right) $ is the identity operator;

\item Outside \del{a}\ins{the} collar neighborhood\del{ of the
boundary}
$X\times \left[0,1\right) $ \ins{of the boundary, $D$} \del{ the
operators under consideration are }\ins{is a}
pseudodifferential operator\del{s} of order $m$;

\item In the domain $X\times \left(\varepsilon ,1\right)
$\ins{,}
\del{the operators are}\ins{$D$ is a}
pseudodifferential operator\del{s} with continuous
symbol\del{s}\footnote{\del{The necessity}\ins{We are forced}
to consider operators with continuous symbols\ins{,} \del{stems from the
observation that}\ins{since}  in the case of pseudodifferential coefficients
$D_k\left( t\right) $ ¢ (\ref{omeg})\ins{,} the\del{
corresponding} symbol \ins{of $D$} is not
\del{in general a }smooth \del{one}\ins{in general}.}
(see \cite{AtSi1} or \cite{ReSc1,Hor3}).%
\end{enumerate}

For this class of operators\ins{,} boundary value problems
\del{are}\ins{can be}  posed in the
same way as\del{ in the} above. The definition of \ins{the} subbundles
$L_{\pm}\left( D\right) $ and the ellipticity condition remain
\del{true}\ins{valid}.

\subsection{Example}
\label{eins}
Let us consider an example of \ins{a} boundary value
problem\del{s} for operators of
the form (\ref{omeg}).

{\em
{\em On a manifold }$M$\ins{,} {\em \ins{we} consider a bundle }$E${
\em
and a decomposition of this bundle in a neighborhood of the boundary
}$X$ {\em into \del{a}\ins{the}  sum of two subbundles}
\begin{equation}
\left. E\right| _{U_{X}}=E_{+}\oplus E_{-}.  \label{deco2}
\end{equation}
{\em For the bundles }$\left. E_{\pm }\right|
_{X}$\ins{, }{\em let us \del{fix}\ins{take}
elliptic first-order operators }$%
\Lambda _{\pm }$ {\em with principal symbol\del{s equal to}}
$\left| \xi ^{\prime }\right| .$
{\em We also choose a first-order operator }$\Lambda _M$
{\em on $M$ with
\del{the }principal symbol} $\left| \xi \right| $ {\em
which acts in the bundle }$E${\em .
In accordance with the decomposition (\ref{deco2}), let us consider
the following first-order elliptic operator in a neighborhood
of the boundary\ins{:}}
\begin{equation}
D_{\pm }=\left( -i\frac \partial {\partial t}+i\Lambda _{+}\right) \oplus
\left( i\frac \partial {\partial t}+i\Lambda _{-}\right) :C^\infty \left(
U_{X},E\right) \rightarrow C^\infty \left( U_{X},E\right) .
\label{mu}
\end{equation}
{\em The \del{equality}\ins{relation}\quad }
\[
L_{+}\left( D_{\pm }\right) =0\oplus \pi ^{*}E_{-},\qquad \pi
:S^{*}X\rightarrow X,
\]
{\em shows that the boundary condition}
\begin{equation}
\left. u_{-}\right| _{X}=g\in C^\infty (X,E_{-})\qquad
\text{{\em for} }u=\left( u_{+},u_{-}\right) \in C^\infty \left(
U_{X},E_{+}\oplus E_{-}\right)   \label{bndcond}
\end{equation}
{\em defines an elliptic boundary value problem for
\ins{the} operator
(\ref{mu}). Let us extend }$D_{\pm }$ {\em
\del{inside}\ins{to the interior of}
the manifold. Consider a cut\dell{-}off function $\chi $ on
}$M$, $0\leq \chi \left( t\right) \leq 1, ${\em
that \ins{is} equal\del{s} \ins{to}  }$1$ {\em
for }$0\leq t\leq 1/3$ {\em and is zero \del{when}\ins{for}}
$t\geq 2/3.$ {\em
The \del{required}\ins{desired}
extension of the operator is given by the formula}
\begin{equation}
D_{\pm }=\chi \left( t\right) \left[ \left( -i\frac \partial {\partial%
t}+i\Lambda _{+}\right) \oplus \left( i\frac \partial {\partial t}+i\Lambda
_{-}\right) \right] +\left( 1-\chi \left( t\right) \right) i\Lambda _M.%
\label{primer}
\end{equation}
{\em The boundary value problem for the operator }$D_{\pm }$ {\em
with \ins{the} boundary condition (\ref{bndcond}) is denoted by }$
{\cal D}_{\pm }${\em. It is well-known (\del{see
}e.g.\ins{,} \ins{see}
\cite{Hor3} or \cite{Bout2}) that this boundary value problem has index zero.
This \del{can be proved}\ins{follows}, for example,
\del{by noting}\ins{from the observation}  that the family of boundary
value problems}
\[
{\cal D}_{\pm }+ip
\]
{\em is an elliptic family in the half-plane $%
{\rm Re}\,p>0$ in the sense of Agranovich--Vishik \cite{AgVi1}.
Consequently, it is invertible for \ins{sufficiently}
large values of the parameter
}$p. ${\em The invertibility of \ins{the} family} $
{\cal D}_{\pm }+ip${\em
\ can be shown directly (see \cite{Hor3}).}

{\em \del{In the case, when}\ins{If}  one
of the bundles }$E_{\pm }$ {\em
coincides with the \del{whole}\ins{entire}  }$E,$ {\em
\ins{then the} corresponding operator
}$D_{\pm }$ {\em is denoted by
}$D_{-}$ {\em or }$D_{+}.$ {\em For example, \ins{the} operator }
$D_{+}$ {\em does not contain boundary conditions.}
%\end{example}
}

\section{The$\,$ homotopy$\,$ classification $\,$of$\,$ boundary$\,$ value problems}

\subsection{Classification of operators of order zero}

In the class of elliptic operators on manifolds with
boundary\dell{,} introduced \del{at}\ins{in}
the end of \ins{the} previous section,
operators of order zero play an important
role, since these operators do not require boundary conditions.

\del{An}\ins{The} abelian group of stabl\del{y}\ins{e}
homotop\del{ic}\ins{y} \ins{classes of}  elliptic zero-order operators is
denoted by $\limfunc{Ell}^0\left(M\right)$.

An elliptic operator $D$ of order zero is \del{an operator of
multiplication}\ins{a bundle isomorphism}  in a
neighborhood of the boundary of $M$ (see
(\ref{omeg}))\dell{,}\ins{;} hence, its principal
symbol defines an element of $K$-theory with compact
support\ins{s:}
\[
\left[ \sigma \left( D\right) \right] \in K\left( T^{*}\left( M\backslash
\partial M\right) \right) .
\]
Thus, \del{a}\ins{we have the}  homomorphism\del{ from the group of
stably homotopic elliptic operators}
\del{qqqqqof order zero is defined}

\begin{equation}
\begin{array}{ccc}
\chi :\limfunc{Ell}^0\left( M\right) & \longrightarrow & K\left(%
T^{*}\left( M\backslash \partial M\right) \right) , \\
&  &  \\
 \left[ D\right] & \mapsto & \ \left[ \sigma \left( D\right) \right] .%
\end{array}
\label{sstar2}
\end{equation}
The following theorem gives the homotopy classification of elliptic
operators of order zero.

\begin{theorem}
The map\ins{ping}  \del{$\chi $}\ins{{\em(\ref{sstar2})}}
is an isomorphism of abelian groups.
\end{theorem}

\noindent {\em Proof. }Let us construct the inverse
map\ins{ping}
\[
\chi ^{\prime }:K\left( T^{*}\left( M\backslash \partial M\right) \right)%
\longrightarrow \limfunc{Ell}^0\left( M\right) .
\]
By \del{means}\ins{virtue}  of \del{a}\ins{the}
natural isomorphism\footnote{$B^{*}M$
is \del{a}\ins{the}  unit \ins{co}ball bundle of $M$
(with respect to some Riemannian metric),
\del{while}\ins{and}
$\partial B^{*}M=S^{*}M\cup \left. B^{*}M\right| _{\partial M}$
is its boundary.}
\[
K\left( T^{*}\left( M\backslash \partial M\right) \right) \simeq K\left(
B^{*}M,\partial B^{*}M\right) ,
\]
the group $K\left( T^{*}\left( M\backslash \partial M\right) \right) $
is \del{a}\ins{the} group of stabl\del{y}\ins{e}
homotop\del{ic}\ins{y} \ins{classes of} elliptic symbols on
$M$\del{ which
are} independent of \ins{the} cotangent variables in a neighborhood of
\del{the boundary }$X$. The map\ins{ping}
$\chi ^{\prime }$ is given by the formula
\[
\chi ^{\prime }\left[ \sigma \right] =\left[ \widehat{\sigma }\right] \in %
\limfunc{Ell}^0\left( M\right) ,
\]
where $\widehat{\sigma }$ is an elliptic pseudodifferential operator of
order zero on $M$ with\del{ the} principal symbol $\sigma
$\del{, and} \ins{such that} near the boundary
$\widehat{\sigma }$ is\del{ equal to} a\del{ multiplication
operator} \ins{bundle homomorphism}. It can be shown
that $\chi ^{\prime}$ is \ins{the} inverse \del{to}\ins{of}
$\chi .$ This proves the theorem.

\subsection{Order reduction: from order one to order zero}

In contrast \del{to}\ins{with}  zero-order operators
considered earlier, operators of order
one \ins{in general} require boundary conditions. Nevertheless, the
homotopy classification is the same in both cases.

\begin{definition}
{\em Boundary value problems }$%
{\cal D}_1$ {\em and }${\cal D}_2$ {\em for operators of order
one are\del{ called} \ins{said to be}} stably
homotopic{\em \dell{,}
if for some operators }$%
{\cal D}_{\pm }$ {\em and }${\cal D}_{\pm }^{\prime }$ {\em
(see Example~\ref{eins})
\ins{the} elliptic boundary value problems }
\[
{\cal D}_1\oplus {\cal D}_{\pm }\quad \text{{\em and}\quad }%
{\cal D}_2\oplus {\cal D}_{\pm }^{\prime }
\]
{\em are homotopic.}
\end{definition}

\del{An}\ins{The} abelian group of stabl\del{y}\ins{e}
homotop\del{ic}\ins{y} \ins{classes of}
elliptic boundary value problems for
operators of order one \del{is}\ins{will be}
denoted by $\limfunc{Ell}^1\left( M\right)$.

\begin{theorem}
\label{th2}\ins{The} \del{''}order\ins{-}increasing\del{''}
map\ins{ping}
\begin{equation}
\times D_{+}:\limfunc{Ell}^0\left( M\right) \longrightarrow
\limfunc{Ell}^1\left( M\right)   \label{chto}
\end{equation}
induced by \ins{the} composition with\del{ operator}
$D_{+}$\dell{,}
is an isomorphism of abelian groups.
\end{theorem}

\begin{remark}
{\em In the proof of the theorem\ins{,} \ins{we give} an explicit formula
for the inverse order reduction map\ins{ping}\del{ is given}}
\[
\varepsilon _1=\left( \times D_{+}\right) ^{-1}:\limfunc{Ell}%
^1\left( M\right) \longrightarrow \limfunc{Ell}^0\left(
M\right) .
\]
\end{remark}

\noindent {\em Proof. } Consider a boundary value problem
$\left( D,B\right) $ for a first-order elliptic operator.

\noindent $1)$ First, we construct a homotopy of \ins{the} restriction
of\del{ operator} $D$ to the boundary $X$ together with a homotopy
of the boundary condition $B$\dell{,} such that the boundary value problem
is deformed to \ins{the} model form (\ref{mu}), (\ref{bndcond}).
According to (\ref{omeg}), \ins{the} operator $D$ on the boundary is equal to
\[
D=\gamma \left( \frac \partial {\partial t}+A\right) ,
\]
where $\gamma $ is an isomorphism of \ins{the} vector bundles
$E$ and $F$. \ins{The} \del{E}\ins{e}llipticity of $D$ is equivalent
to the absence of pure\del{ly} imaginary eigenvalues of the principal
symbol of\del{ operator} $A$ for $\left|\xi ^{\prime }\right| =1.$

STEP 1. Let $P$ be a pseudodifferential operator on\del{ the
boundary} $X$
with\del{ the} principal symbol equal to the projection
on\del{to} \ins{the} subbundle
$L_{+}\left( D\right)$ along the complementary \ins{bundle}
$L_{-}\left( D\right)$.
Consider the following homotopy with parameter
$\tau\in\left[0,1\right]$\ins{:}
\begin{equation}
D_\tau ^{\prime }=\gamma \left( \frac \partial {\partial t}+\left( 1-\tau
\right) A+\tau \Lambda _X\left( 2P-1\right) \right). \label{odin}
\end{equation}
\del{During t}\ins{T}his homotopy \ins{takes}
the eigenvalues of the symbol $\sigma \left( A_\tau
\right) $\del{ are deformed} to $\pm 1$ according to \ins{the} formula
\[
\left( 1-\tau \right) \lambda +\tau \limfunc{sign}\lambda ,
\]
while the subspaces $L_{\pm }\left( D_\tau ^{\prime }\right) $
do not change. As a consequence, the operators $D_\tau $
\ins{always} remain elliptic.
\del{During t}\ins{T}he homotopy (\ref{odin}) \ins{does not
change} the boundary symbol
$$
L_{+}\left( D_\tau ^{\prime }\right)
\stackrel{\sigma \left( B\right) }{\rightarrow }\pi
^{*}G.
$$

STEP 2. Let us embed the bundle $G$ of boundary values in a trivial bundle
${\bf C}^N\supset G$ and denote by $P_G$ the projection
on\del{to}
sections of $G$ in the space $C^\infty \left(
X,E\oplus {\bf C}^N\right)$.
Consider \del{a}\ins{the
following}  homotopy of almost-projections
with parameter $\varphi \in \left[ 0,\pi/2\right] $\ins{:}
\begin{equation}
P_\varphi =P\limfunc{cos}^2\varphi +P_G\sin ^2\varphi +\limfunc{cos}%
\varphi \sin \varphi \left( BP+B^{-1}P_G\right) .  \label{mnogodot}
\end{equation}
The almost-projections act in the
direct sum $C^\infty \left(X,E\oplus{\bf C}^N
\right)$.
In this homotopy\ins{,} the bundle defined by the projection
$\sigma \left( P_\varphi \right) $ is \del{a result of}\ins{the}  rotation
by \ins{the} angle $\varphi $ of \ins{the}
subbundle $L_{+}\left( D\right) \subset
\pi ^{*}\left(E\oplus {\bf C}^N\right) $ towards
the subbundle $\func{Im}\sigma \left( P_G\right) =\pi ^{*}G$ with the help of
the isomorphism $\sigma \left( B\right) .$ This homotopy
defines \del{a}\ins{the}  homotopy of
operators
\begin{equation}
D_\varphi ^{\prime }=\gamma \left( \frac \partial {\partial t}+\Lambda
_X\left( 2P_\varphi -1\right) \right)   \label{dva}
\end{equation}
and \del{a}\ins{the}   homotopy of boundary conditions
\begin{equation}
\begin{array}{c}
B_\varphi :C^\infty \left( X,\left. E\right| _X\oplus {\bf C}^N\right)
\longrightarrow C^\infty \left( X,G\right) \subset C^\infty \left( X,\left.
E\right| _X\oplus {\bf C}^N\right) , \\
\\
B_\varphi =\limfunc{cos}\varphi BP+\sin \varphi P_G.
\end{array}
\label{xxstar}
\end{equation}
For the final value of the parameter\ins{,} $\varphi =\pi
/2$\ins{,} we obtain
\[
\begin{array}{c}
D_{\pi /2}^{\prime }=\gamma \left( \frac \partial {\partial t}+\Lambda
_X\left( 2P_G-1\right) \right) , \\
\\
B_{\pi /2}=P_G:C^\infty \left( X,\left. E\right| _X\oplus {\bf C}%
^N\right) \longrightarrow C^\infty \left( X,G\right) ,
\end{array}
\]
which coincides with the model operator ${\cal D}_{\pm }$
up to an isomorphism of vector bundles.

\vspace{0.2cm}

\noindent 2) T\ins{he t}wo homotopies (\ref{odin})
and (\ref{dva}), (\ref{xxstar})
of the restriction of \ins{the} boundary value problem to the boundary
can be lifted\del{ up} to a homotopy of boundary value problems.
To this end, \ins{we} consider a cut\del{-}off function
$\psi :M\rightarrow {\bf R}$ that \ins{is} equal\del{s}
\ins{to} one
in a neighborhood of the boundary of $M$\dell{,} and is\del{
equal to} zero outside
the domain $X\times \left[ 0,1/2\right) $.
The composition of \ins{the} homotopies (\ref{odin})\dell{,}
\ins{and} (\ref{dva}) is denoted
for \del{short}\ins{brevity}  by
$\left( D_\tau ^{\prime },B_\tau \right)$,
$\tau \in \left[ 0,1\right]$.
Let us attach a finite cylinder to the manifold $M$
(see Fig.~\ref{ff1})\ins{:}

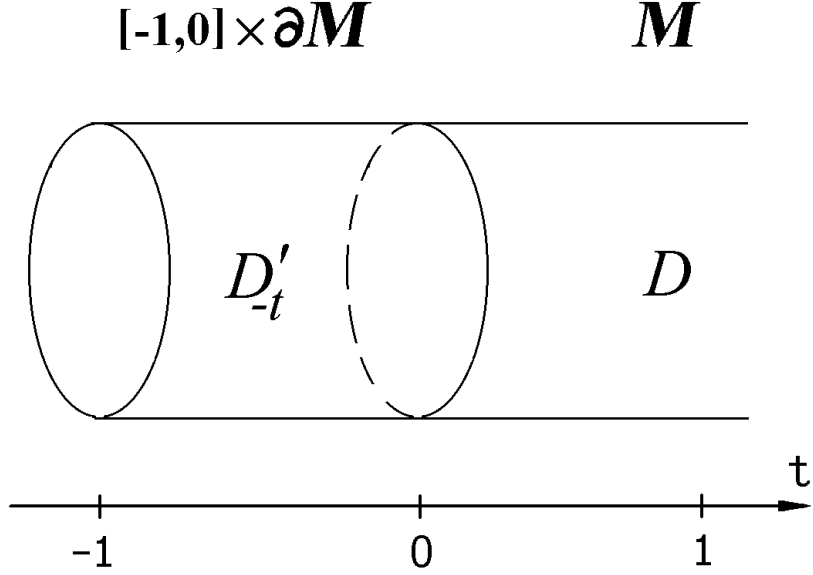
\fgr{fig.eps}{ff1}{The operator  $D$  on the manifold
$M'=(\left[ -1,0\right]\times \partial M) \cup M$}

\[
M^{\prime }=\left[ -1,0\right]\times X\cup M. %
\]
\ins{The} \del{O}\ins{o}perator $D$ \ins{can be}
extend\del{s}\ins{ed}  to this manifold: on the cylinder
$\left[ -1,0\right]\times X $\ins{,} it is defined by the homotopy
$D_{-t}^{\prime }$.
The required lifting of the homotopy
$\left( D_\tau ^{\prime},B_\tau \right)$
\del{up }to a homotopy of boundary value problems\
$\left( D_\tau,B_\tau \right) $ on $M$ is defined by the formula
\begin{equation}
D_\tau \left( t,x\right) =D\left( t-\psi \left( t\right) \tau ,x\right)
,\qquad \tau \in \left[ 0,1\right] .  \label{kapp}
\end{equation}

Thus, the boundary value problem $\left(D,B\right) $ is now deformed to
a boundary value problem $\left(D_1,B_1\right)$\dell{,}
\del{which}\ins{that}  coincides near the
boundary with the model problem ${\cal D}_{\pm }.$
Hence, \ins{we have defined} \del{an}\ins{the} \ins{zero-order}
elliptic operator\del{ of order zero is defined}
\[
\left[ {\cal D}_{\pm }^{-1}\circ \left( D_1,B_1\right) \right] \in %
\limfunc{Ell}^0\left( M\right) .
\]
It can be verified that this construction defines a homomorphism of groups
\[
\begin{array}{ccc}
\varepsilon _1:\limfunc{Ell}^1\left( M\right) & \longrightarrow &
\limfunc{Ell}{}^0\left( M\right), \\
\left[ D,B\right] & \longmapsto & \left[ {\cal D}_{\pm }^{-1}\circ \left(%
D_1,B_1\right) \right].
\end{array}
\]
Indeed, this construction is uniquely determined; it takes
direct sums of boundary value problems to sums of the
corresponding elements; the model operators are taken to zero;
finally, the construction is homotopy invariant.
\ins{It follows} \del{F}\ins{f}rom the definition of stable
homotop\del{y}\ins{ies}  for boundary value problems
\del{it follows }that this map\ins{ping}
is the inverse \del{to}\ins{of} the order\ins{-}increasing map\ins{ping}
(\ref{chto}). This establishes the reduction of classical
boundary value problems of order one to operators of order zero.
The theorem is \ins{thereby} proved.

\subsection{Order reduction: from an arbitrary order to order one}

\begin{definition}
{\em Elliptic boundary value problems }${\cal D}_1$ {\em and }
${\cal D}_2$ {\em of order }$m\geq 2 ${\em
are \del{called}\ins{said to be} }
stably homotopic {\em if for some model operators\ }
${\cal D}_{\pm }$ {\em and }${\cal D}_{\pm }^{\prime }$ {\em
there exists a homotopy between \ins{the} boundary value problems}
\[
{\cal D}_1\oplus {\cal D}_{\pm }D_{+}^{m-1}\quad \text{{\em and}}%
\qquad {\cal D}_2\oplus {\cal D}_{\pm }^{\prime }D_{+}^{m-1}.
\]
\end{definition}

The group of stabl\del{y}\ins{e} homotop\del{ic}\ins{y}
\ins{classes of}  boundary value problems
for operators of order $m$ is
denoted by $\limfunc{Ell}^m\left( M\right)$.

\begin{theorem}
\label{th3}The map\ins{ping}\del{ of order increase by $m-1$}
\[
\left( \times D_{+}\right) ^{m-1}:\limfunc{Ell}^1\left( M\right)
\longrightarrow \limfunc{Ell}^m\left( M\right),
\]
\ins{which increases the order by $m-1$,}
is an isomorphism of \del{A}\ins{a}belian groups.
\end{theorem}

\noindent {\em Proof. }Consider a boundary value problem
$\left( D,B\right) $ for an operator $D$ of order $m.$
The direct sum
\begin{equation}
\left( D,B\right) \oplus \bigoplus\limits_1^{m-1}D_{+}^m  \label{aux1}
\end{equation}
defines the same element in the group $\limfunc{Ell}^m\left( M\right)$
as the original problem $\left( D,B\right)$.
Let us construct a homotopy of \ins{the} boundary value problem
(\ref{aux1}) to \del{a}\ins{the}
composition of a boundary
value problem for an operator of order one and \ins{the} operator
$D_{+}^{m-1}.$ As in the proof of Theorem \ref{th2} (see (\ref{kapp})),
it suffices to construct \del{a}\ins{the corresponding}
homotopy of the restriction of\del{ the
operator} $D$ to the boundary together with the
boundary conditions\ins{.}\del{ which deforms them to the required
form.}

Let us represent\del{ operator} $D$ in the form
\[
D=\sum\limits_{k=0}^mD_k\left( -i\frac \partial {\partial t}-i\Lambda
_X\right) ^k\left( -i\frac \partial {\partial t}+i\Lambda _X\right)
^{m-k}\equiv \sum\limits_{k=0}^mD_kD_{-}^kD_{+}^{m-k},
\]
where $\Lambda _X $ is\del{, as above,} \ins{again}
a first-order operator with principal
symbol\del{ equal to} $\left| \xi ^{\prime }\right|$\dell{,}
and \ins{the} $D_k$ are \ins{zero-order}
pseudodifferential\del{ zero-order} operators on $X$.
By virtue of condition (\ref{omeg}), \del{it}\ins{we} can\del{
be} assume\del{d} that the sum
\[
\sum\limits_{k=0}^mD_k
\]
is equal to the identity operator. Consider
\del{an}\ins{the}  operator homotopy
\[
D_\tau ^{\prime }=\left(
\begin{array}{cccc}
D+\tau ^mD_0D_{+}^m-\tau ^mD & \tau ^{m-1}D_1D_{+}^m & \ldots & \tau
D_{m-1}D_{+}^m+\tau D_mD_{-}D_{+}^{m-1} \\
-\tau D_{-}D_{+}^{m-1} & D_{+}^m & \ldots & 0 \\
0 & -\tau D_{-}D_{+}^{m-1} & \ldots & 0 \\
\ldots & \ldots & \ddots & \ldots \\
0 & 0 & \ldots & D_{+}^m
\end{array}
\right) .
\]
At the initial point $\tau =0$\ins{,} we have
\[
D_{\tau =0}^{\prime }=D\oplus \bigoplus\limits_1^{m-1}D_{+}^m.
\]
On the other hand, \del{for}\ins{at}  $\tau =1$ \ins{we obtain}
the factorization required in
the theorem\ins{:}\del{ is obtained}
\[
D_{\tau =1}^{\prime }=\left(
\begin{array}{cccc}
D_0D_{+} & D_1D_{+} & \ldots & D_{m-1}D_{+}+D_mD_{-} \\
-D_{-} & D_{+} & 0 & 0 \\
0 & -D_{-} & D_{+} & \ldots \\
0 & 0 & \ldots & D_{+}
\end{array}
\right) \circ mD_{+}^{m-1}. %
\]
The coefficient of $\left( -i\frac \partial {\partial t}\right) ^m$ in
the operator $D_\tau ^{\prime }$ is equal to \del{a}\ins{the} composition
\[
D_{\tau m}^{\prime }=\left(
\begin{array}{cccc}
1 & \tau ^{m-1}\left( D_1+\ldots +D_m\right) & \ldots & \tau \left(
D_{m-1}+D_m\right) \\
0 & 1 & \ldots & 0 \\
0 & 0 & 1 & \ldots \\
0 & 0 & \ldots & 1
\end{array}
\right) \circ \left( %
\begin{array}{cccc}
1 & 0 & \ldots & 0 \\
-\tau & 1 & \ldots & 0 \\
0 & -\tau & 1 & \ldots \\
0 & 0 & \ldots & 1
\end{array}
\right) .
\]
Thus, for the operators
\[
D_\tau =\left( D_{\tau m}^{\prime }\right) ^{-1}D_\tau ^{\prime }
\]
the corresponding coefficient is equal to
\del{identity}\ins{unity}. Let us show that
the operator $D_\tau $ is elliptic for $\tau \in \left[ 0,1\right] $.
To this end, we compute the subspace $L_{+}\left( D_\tau
\right)$.

Consider a bounded solution
\[
U=\left( U_0\left( t\right) ,U_1\left( t\right) ,\ldots ,U_{m-1}\left(
t\right) \right),\qquad t\rightarrow +\infty,
\]
of the equation
\begin{equation}
\sigma \left( D_\tau \right) \left( x,\xi ^{\prime },-i\frac d{dt}\right)
U=0,\qquad \text{ for }\left| \xi ^{\prime }\right| =1.\text{ }
\label{tristar}
\end{equation}
\del{The e}\ins{E}quation~(\ref{tristar}) can be replaced by an equivalent equation
with the symbol of the operator $D_\tau ^{\prime }$.
\ins{The} \del{B}\ins{b}ounded function $U$
is a solution of an ordinary differential equation
with constant coefficients\dell{,}\ins{;} hence, its derivatives are
\ins{also} bounded\del{ as well}.
\ins{The} \del{C}\ins{c}omponent\dell{-}wise
representation of (\ref{tristar}) \del{leads to}\ins{gives} the system
\begin{equation}
\left\{
\begin{array}{l}
\left( 1-\tau ^m\right) \sigma \left( D\right) \left( -i\frac d{dt}\right)
U_0+\left( -i\frac d{dt}+i\right) ^m\left( d_0\tau ^mU_0+\ldots +d_{m-1}\tau
U_{m-1}\right) + \\
+\tau \left( -i\frac d{dt}+i\right) ^{m-1}\left( -i\frac d{dt}-i\right)
d_mU_{m-1}=0, \\
\left( -i\frac d{dt}+i\right) ^mU_j=\tau \left( -i\frac d{dt}+i\right)
^{m-1}\left( -i\frac d{dt}-i\right) U_{j-1},\qquad \text{ }%
0<j<m
\end{array}
\right.  \label{krut}
\end{equation}
(\ins{the} $d_j$ \del{denote}\ins{are the}
principal symbols of \ins{the} operators $D_j$).
The equation
\[
\left( -i\frac d{dt}+i\right) u=0
\]
on the half-line $\left\{ t\geq 0\right\} $ has only a trivial bounded
solution. Hence, the operator $-id/dt+i$  can be cancel\del{l}ed in
(\ref{krut}) in all equations\dell{,} except \ins{for} the
first\del{ one}. Consequently,
\[
\left( -i\frac d{dt}+i\right) U_j=\tau \left( -i\frac d{dt}-i\right)
U_{j-1}.
\]
Substituting these relations\del{ one} into \ins{one} another, we obtain
\[
\left( -i\frac d{dt}+i\right) ^jU_j=\tau ^j\left( -i\frac d{dt}-i\right)
^jU_0.
\]
\del{This implies}\ins{It follows} that the first
equation in (\ref{krut}) is reduced
to the requirement
\begin{equation}
\sigma \left( D\right) \left( -i\frac d{dt}\right) U_0=0.  \label{reduk}
\end{equation}
\del{Thus,}\ins{We conclude that} the operator $D_\tau $ is
indeed\del{ an} elliptic\del{ one}, since
equation (\ref{reduk}) has no \ins{solutions} bounded on the
\del{whole}\ins{entire}  line\del{ solutions}.

Hence, \ins{we have obtained}
the following description of the bundle $L_{+}\left(
D_\tau\right) $\ins{:}
\del{is obtained: }the projection on the first term
\del{of}\ins{in}  the sum
\[
E\oplus \bigoplus\limits_1^{m-1}E\stackrel{pr}{\longrightarrow }E
\]
induces an isomorphism of vector bundles
\[
L_{+}\left( D_\tau \right) \stackrel{pr}{\longrightarrow }L_{+}\left(
D\right) ;
\]
the preimage of an element $u\in L_{+}\left( D\right) $
under this map\ins{ping} is given by the formula
\begin{eqnarray}
U &=&\left( U_0,\ldots ,U_{m-1}\right) ,  \label{tozhd} \\
U_0 &=&u,  \nonumber \\
\left( -i\frac d{dt}+i\right) ^jU_j &=&\tau ^j\left( -i\frac d{dt}-i\right)
^jU_0.  \nonumber
\end{eqnarray}
Let us decompose the operator of boundary conditions
in the same way\dell{,} as\del{ the operator} $D$\del{
has been rewritten in the above}:
\[
Bj_X^{m-1}=\left. \sum\limits_{k=0}^{m-1}B_k\left( -i\frac \partial
{\partial t}-i\Lambda _X\right) ^k\left( -i\frac \partial {\partial
t}+i\Lambda _X\right) ^{m-1-k}\right| _{t=0}.
\]
This implies that the boundary condition
\[
Bj_X^{m-1}\circ pr:C^\infty \left( M,E\oplus%
\bigoplus\limits_1^{m-1}E\right) \rightarrow C^\infty \left( X,G\right)
\]
has a factorization: on the subspace
$L_{+}\left( D_{\tau =1}\right) $, by virtue of (\ref{tozhd}), we have
\[
\sigma \left( B\right) j_X^{m-1}\circ pr=\sigma \left( B^{\prime }\right)%
j_X\circ \left( -i\frac d{dt}+i\right) ^{m-1},\qquad \text{where }%%
\left( \sigma \left( B^{\prime }\right) j_X\right)
U=\sum\limits_{k=0}^{m-1}b_kU_k\left( 0\right) .
\]
That is why the homotopy of boundary value problems
\begin{equation}
\left( D_\tau ,B\circ pr\right) ,\qquad \tau \in \left[ 0,1\right],
\label{odinx}
\end{equation}
connects the initial problem (\ref{aux1}) with the composition of a boundary
value problem for a first-order operator and \ins{the}
operator $D_{+}^{m-1}$\ins{:}
\begin{equation}
\left( D_{\tau =1},B\circ pr\right) =\left( D^{\prime },B^{\prime }\right)%
\circ mD_{+}^{m-1}.  \label{xx}%
\end{equation}

\del{It}\ins{One} can\del{ be} show\del{n} that
the correspondence \del{of}\ins{between the}  boundary value problems
$\left(D,B\right) $ and $\left( D^{\prime },B^{\prime }\right) $
for operators of order one induces a map\ins{ping}
\[
\varepsilon _m^{\prime }:\limfunc{Ell}^m\left( M\right)
\longrightarrow \limfunc{Ell}^1\left( M\right) .
\]
Let us check that this map\ins{ping} is the inverse
\del{to}\ins{of}  the map\ins{ping}
$\times D_{+}^{m-1}.$ Indeed, the homotopy (\ref{odinx}) shows
that the group $\limfunc{Ell}^m\left( M\right) $ is generated by compositions
of operators (\ref{xx}), i.e. by the range of the
map\ins{ping} $\times D_{+}^{m-1}.$
Hence, this map\ins{ping} is onto\del{ map}. Let us prove
\del{the equality}\ins{that}
\[
\varepsilon _m^{\prime }\circ \left( \times D_{+}^{m-1}\right) =Id.%
\]
Indeed, for an elliptic boundary value problem
$\left( D^{\prime },B^{\prime }\right) \circ D_{+}^{m-1}$ of
order $m$\ins{,}
the matrix of the operator $D_\tau $ in the homotopy
(\ref{odinx}) is a product of two triangular matrices with constant
\del{coefficients on the }diagonal \ins{entries}
(with respect to the parameter $\tau $
of the homotopy). Thus, this homotopy is trivial, i.e.\ins{,} homotopic to a
constant homotopy.

Theorem \ref{th3} is \ins{thereby} proved.

\subsection{Main theorems}

The \ins{above} results on the homotopy
classification of boundary value problems of
fixed order\del{s obtained in the above} are summarized in the following
theorems.

\begin{theorem}
\label{th4}{\em (\ins{the} Atiyah--Bott
obstruction to the existence of elliptic
boundary value problems)}
For an elliptic operator $D$ on a manifold $M$ with boundary
$X$,
the following conditions are equivalent.
\end{theorem}

\begin{enumerate}
\item  {\em The operator }$D$
{\em stably,
i.e. up to the direct sum with an operator of the form
 $D_+^{m-1}{\cal D}_\pm$ (cf.~Definitions 1, 2),
admits an elliptic boundary value problem;}

\item {\em The following inclusion holds\ins{:}}
\[
\left[ L_{+}\left( D\right) \right] \in \pi ^{*}K\left( X\right),
\quad\pi
:S^{*}X\rightarrow X;
\]

\item  {\em The restriction }$\!
\left. \sigma \left( D\right) \right| _X$ {\em of the principal symbol
of the operator to the boundary }$\!X\!$
{\em is stably homotopic to \del{a}\ins{the}
symbol of \ins{a} multiplication operator\/};

\item  $j^{*}\left[ \sigma \left( D\right) \right] =0$ {\em for%
} $j^{*}:K\left( T^{*}M\right)\! \rightarrow\! K\left( T^{*}X\times
{\bf R}%
\right) =K^1\left( T^{*}X\right) ,$ {\em where}
$j:\left.T^*M\right|_X\rightarrow T^*M$
{\em is the inclusion}.
\end{enumerate}

\noindent {\em Proof}.
The equivalence of conditions 1) and 2) follows from the definition
of ellipticity for boundary value problems.
The equivalence of 3) and 4) is a consequence of the definition
of the group $K\left( T^{*}X\times {\bf R}\right) $ in terms
of the difference construction.

Let us check the equivalence of conditions 2) and 4).
By virtue of homotopies constructed in Theorems \ref{th2} and \ref{th3},
it can be assumed that the operator $D$ in a neighborhood of
the boundary has the form
\begin{equation}
D=\frac \partial {\partial t}+\Lambda _X\left( 2P-1\right) .
\label{doublestar}
\end{equation}
For the operator (\ref{doublestar})\ins{,} the following
formula is valid\ins{:}
\[
j^{*}\left[ \sigma \left( D\right) \right] =\delta \left[ L_{+}\left(
D\right) \right] ,\quad \delta :K\left( S^{*}X\right) \rightarrow K^1\left(
T^{*}X\right).
\]
The kernel of \ins{the} homomorphism $\delta $ \del{consists
exactly of}\ins{coincides with}  the
subgroup $\pi ^{*}K\left( X\right)
\subset K\left( S^{*}X\right) $ \cite{APS3}.
This implies the equivalence of 2) and 4).
\ins{The} \del{T}\ins{t}heorem is \ins{thereby} proved.

\begin{theorem}
{\em (\ins{the} homotopy classification of elliptic boundary value problems)}
\ins{For $m\geq1$,} \del{T}\ins{t}here is
an isomorphism of groups\del{ for $m\geq 1$}
\[
\varepsilon _m:\limfunc{Ell}^m\left( M\right) \longrightarrow
\limfunc{Ell}^0\left( M\right)
\]
\del{which}\ins{that}  is \ins{the} inverse \del{to}\ins{of the}
order\ins{-}increasing map\ins{ping}
\[
\times D_{+}^m:\limfunc{Ell}^0\left( M\right) \longrightarrow
\limfunc{Ell}^m\left( M\right) .
\]
Moreover, the \ins{following} symbol isomorphism holds\ins{:}
\[
\chi :\limfunc{Ell}^0\left( M\right) \longrightarrow K\left(%
T^{*}\left( M\backslash \partial M\right) \right) .
\]
\end{theorem}

\begin{corollary}
For elliptic boundary value problems ${\cal D}$, \del{there is
an}\ins{one has the}  \del{equality}\ins{equation}
\[
\limfunc{ind}{\cal D}=\limfunc{ind}\varepsilon _m\left[ {\cal D}%
\right] ,\qquad m\geq 1,
\]
and \del{an}\ins{the}  index formula
\begin{equation}
\limfunc{ind}D=p_{!}\left[ \sigma \left( D\right) \right] ,\quad
\label{exci}
\end{equation}
where  $$p:M\to pt$$ and
\[
\left[ \sigma \left( D\right) \right] \in K\left( T^{*}\left( M\backslash
\partial M\right) \right)
\]
for \ins{an} operator $D$ of order zero representing
$\varepsilon_m[{\cal D}]$.
\end{corollary}

\noindent {\em Proof. }The map\ins{ping} $\varepsilon _m$ preserves
\ins{the} index by definition. \del{The
equality}\ins{Equation}~(\ref{exci}) is a special case of the
``excision" property  of the index (see \cite{AtSi1}).

\begin{corollary}
\label{sled1}{\em (cobordism invariance of \ins{the} index)}
Let $M$ be a compact manifold with boundary $X$. \ins{We}
\del{D}\ins{d}enote the natural
inclusion map\ins{ping} by\del{ $j$}
$$
j:X \longrightarrow M.
$$
Consider the induced map\ins{ping}
$$
j^{*}:K^1\left( T^{*}M\right) \rightarrow K\left(T^{*}X\right)
$$
\ins{in  $K$-theory.}
If an elliptic operator $D$ over\del{ the boundary}
$X$ satisfies the inclusion
\[
\left[ \sigma \left( D\right) \right] \in \func{Im}j^{*},
\]
then
\[
\limfunc{ind}D=0.
\]
\end{corollary}

\noindent {\em Proof. }The \del{requ}\ins{des}ired
statement follows from the commutative
diagram
\[
\begin{array}{ccccc}
K^1\left( T^{*}M\right) & \rightarrow & K\left( T^{*}X\right) & \rightarrow
& K\left( T^{*}\left( M\backslash \partial M\right) \right) \\
&  & \chi \uparrow \; &  & \uparrow \\ %
&  & \limfunc{Ell}\left( X\right) & \stackrel{\alpha }{\rightarrow } &
\limfunc{Ell}^1\left( M\right) .
\end{array}
\]
\ins{Here} \del{T}\ins{t}he upper row\del{ here}
is induced by the exact sequence of \ins{the} triple
\[
S^{*}M\subset \partial B^{*}M\subset B^{*}M,
\]
$\chi $ \del{denotes}\ins{stands for the}
difference construction on \ins{the} (closed) manifold
$X$, and the map\ins{ping} $\alpha $ \del{associates to
an}\ins{takes each}  elliptic operator
\[
B:C^\infty \left( X,{\bf C}^N\right) \rightarrow C^\infty \left(
X,G\right)
\]
\del{a}\ins{to the}  boundary value problem
\[
\left\{
\begin{array}{ll}
D_{-}u=f, & \quad u,f\in C^\infty \left( M,{\bf C}^N\right) ,
\vspace{0.2cm} \\
Bj_Xu=g, & \ \quad g\in C^\infty \left( X,G\right) .
\end{array}
\right.
\]

\section{Boundary value problems for general elliptic equations}

\subsection{Spectral boundary value problems}
For an arbitrary elliptic operator $D$, which in general does not satisfy
\ins{the} Atiyah--Bott condition (see Theorem \ref{th4}
\del{of}\ins{in}  the previous section),\del{
in \cite{ScSS18}
there were introduced}
boundary value problems of
the following form
\ins{were introduced in \cite{ScSS18}:}
\begin{equation}
\left\{
\begin{array}{ll}
Du=f, & u\in H^s\left( M,E\right), \quad f\in H^{s-m}\left( M,F\right) ,
\vspace{0.2cm}  \\
Bj_X^{m-1}u=g,\qquad & g\in \func{Im}P\subset H^\sigma \left( X,G\right) ,
\end{array}
\right.  \label{star1a}
\end{equation}
where \ins{the subspace} $\func{Im}P$ is
\del{a subspace of Sobolev space on the boundary,
defined by}\ins{the range of}  a pseudodifferential projection $P$
of order zero \ins{in a Sobolev space on the boundary}.
It \del{is}\ins{was} also shown in \cite{ScSS18}
that the boundary value problem
(\ref{star1a}) \del{has}\ins{is}  Fredholm\del{ property}
if and only if it is {\em elliptic},
i.e.\ins{,} its {\em boundary symbol}
\[
\sigma \left( B\right) :L_{+}\left( D\right) \rightarrow \func{Im}\sigma
\left( P\right)
\]
is a bundle isomorphism. This class of boundary value problems does not
carry obstructions of \ins{the} Atiyah--Bott type,
since for an arbitrary elliptic
operator $D$ there exists \del{the}\ins{a}
so-called {\em spectral boundary value problem},
which has the Fredholm property \cite{NScSS3}.

\noindent{\bf Example }
{\em \ins{The} \del{S}\ins{s}pectral boundary value
problem for an operator of order one.}

Let  $D$ be a first-order elliptic operator. In a neighborhood
of the boundary\ins{,} it has the form
\[
D=\gamma \left( t\right) \left( \frac \partial {\partial
t}+A\right),
\]
where $\gamma(t)$ is a bundle isomorphism.
The ellipticity of $D$ implies that the principal symbol
$\sigma \left( A\right) \left( x,\xi ^{\prime }\right) $ has no
pure\del{ly} imaginary eigenvalues for $\left| \xi^{\prime }\right| =1.$
Thus, the family
\[
ip+A
\]
is\del{ an} elliptic\del{ family}
in the sense of Agranovich--Vishik in \del{a}\ins{some}
sector\del{ surrounding} \ins{containing}
the real line $\left\{ p\in {\bf R}\right\} $.
It \del{is}\ins{was}  proved in \cite{NScSS3}
that the spectral projection $P_{+}$
of the operator $A$ on\del{to} the subspace corresponding to spectral points
with nonnegative real parts along the corresponding negative subspace
is a pseudodifferential projection. Its principal symbol
$\sigma \left( P_{+}\right) $ is equal to the nonnegative spectral projection
for the principal symbol of\del{ the operator} $A$\ins{:}
\begin{equation}
\sigma \left( P_{+}\left( A\right) \right) =P_{+}\sigma \left( A\right) .
\label{mu4}
\end{equation}

\begin{definition}
\ins{The} \del{S}\ins{s}pectral boundary value problem
{\em (cf. \cite{APS1}) for \ins{the} operator}
$D$ {\em is \del{a}\ins{the}  system of equations
of the form}
\begin{equation}
\left\{
\begin{array}{ll}
Du=f, & u\in H^s\left( M,E\right) ,f\in H^{s-1}\left( M,F\right) ,
\vspace{0.2cm} \\
P_{+}\left. u\right| _X=g, & g\in \func{Im}P_{+}.
\end{array}
\right.   \label{lastx}
\end{equation}
\end{definition}

This boundary value problem has \ins{the} Fredholm property, since
\del{the equality}\ins{Eq.}~(\ref{mu4}) implies
that its boundary symbol is\del{ equal to} the identity
\ins{mapping}
\[
\func{Im}P_{+}\sigma \left( A\right) \stackrel{Id}{\longrightarrow }\func{Im}%
\sigma \left( P_{+}\right) ,\qquad L_{\pm }\left( D\right) =\func{Im}P_{\pm
}\sigma \left( A\right).
\]

\subsection{The reduction theorem}
The group of stably homotopic boundary value problems (\ref{star1a}) for
operators of order $m\geq 1$ \del{is}\ins{will be}  denoted by
$\limfunc{Ell}^m\left( M,\partial M\right)$, and the group of stably
homotopic spectral boundary value problems for first-order operators
\del{is}\ins{will be}
denoted by $\limfunc{Spec}\left( M,\partial M\right) .$
Here homotopies are families of boundary value problems
(\ref{star1a}) such that $D$, $B$, and $P$ continuously depend
on the parameter and the trivial problems used in
stabilization are the same as in the case of classical
boundary value problems, i.e., have the form ${\cal D}_\pm
D_+^{m-1}$.

The violation of \ins{the} Atiyah--Bott condition
makes it impossible to reduce boundary
value problems to zero-order operators. Nevertheless,
the homotopies of \ins{the} classical theory, described in
\del{s}\ins{S}ection 2,
can be generalized to the present situation. They result in the following
theorems.

\begin{theorem}
\label{th6}\ins{A} \del{B}\ins{b}oundary value problem
for an operator of order $m\geq 2$ \ins{can be}
reduce\del{s}\ins{d}
to a first-order boundary value problem\dell{,}\ins{.}
\del{i.e.}\ins{In other words,}  there is an isomorphism of
groups
\[
\limfunc{Ell}^m\left( M,\partial M\right) \longrightarrow \limfunc{%
Ell}^1\left( M,\partial M\right)
\]
that is \ins{the} inverse \del{to}\ins{of the}
order\ins{-}increasing map\ins{ping}
\[
\times D_{+}^{m-1}:\limfunc{Ell}^1\left( M,\partial M\right)
\longrightarrow \limfunc{Ell}^m\left( M,\partial M\right) .
\]
\end{theorem}

\begin{theorem}
\label{th7}\ins{A} \del{B}\ins{b}oundary value problem for a
first\ins{-}order operator can be
reduced to a spectral boundary value problem\dell{,}\ins{.}
\del{i.e.}\ins{In other words,}  there is an isomorphism
of groups
\[
\limfunc{Ell}^1\left( M,\partial M\right) \longrightarrow \limfunc{%
Spec}\left( M,\partial M\right) .
\]
\end{theorem}

The proof of Theorem \ref{th6} coincides
with \del{the proof}\ins{that} of \ins{the} similar theorem
(Theorem~\ref{th3}) for
classical boundary value problems, since the formulas given there do not
take into account the classical type of boundary value problems.
\vspace{0.2cm}

\noindent {\em Proof\/} of Theorem \ref{th7}. Consider
\ins{the} boundary
value problem (\ref{star1a}). The first homotopy  (\ref{odin})
\del{of}\ins{in}  the proof of Theorem \ref{th2}
\ins{can be} generalize\del{s}\ins{d}  without
change\ins{s}. Let us substitute \ins{the} pseudodifferential projection $P$ that
defines the boundary values in\ins{to}
the rotation homotopy (\ref{mnogodot})
instead of\del{ projection} $P_G$. \del{As a result}\ins{In
the end}  of \ins{the} homotopy
(\ref{dva}), (\ref{xxstar})\ins{,} we obtain \ins{the}
spectral boundary value problem
\[
\begin{array}{c}
D_{\pi /2}^{\prime }=\gamma \left( \frac \partial {\partial t}+\Lambda
_X\left( 2P-1\right) \right) , \\
\\
B_{\pi /2}=P:C^\infty \left( X,\left. E\right| _X\oplus {\bf C}^N\right)
\longrightarrow \func{Im}P\subset C^\infty \left( X,\left. E\right| _X\oplus
{\bf C}^N\right) .
\end{array}
\]
Thus, an arbitrary first-order boundary value problem
\ins{can be} reduce\del{s}\ins{d}  to
a spectral boundary value problem\dell{,} \del{such that
the}\ins{whose}  spectral subspace coincides
with the subspace of boundary values of the initial problem. This proves
the theorem.

\vspace{0.1cm}

In the general case\ins{,}
\del{a}\ins{the}  reduction of \del{the}\ins{a}
boundary value problem to an operator
of order zero is impossible by \ins{the} Atiyah--Bott condition. In the next
section\ins{,} we discuss a class of boundary value
problems\del{, where} \ins{for which the}
Atiyah--Bott condition is satisfied rationally. \ins{The}
\del{R}\ins{r}eduction to\del{ a} classical
boundary value problem\ins{s}  is \del{done}\ins{carried out}
(also rationally) in this case.

\section{Boundary value problems in even and odd subspaces}

\subsection{Parity conditions}

On the cotangent bundle of the manifold $M$\ins{,} \ins{we}
consider the antipodal
involution\del{ $\alpha $}
$$
\alpha :T^{*}M\longrightarrow T^{*}M,\qquad
\\
 \alpha \left( x,\xi \right)=\left( x,-\xi \right) .
$$

\begin{definition}
{\em A pseudodifferential projection }$P$ {\em of order zero
is \del{called}\ins{said to be}  }even\/ {\em(}odd\/{\em)}{\em \dell{,}
if its homogeneous principal symbol
on the sphere bundle }$S^{*}M$ {\em is invariant (antiinvariant)
with respect to the involution }$\alpha $:
$$
\alpha ^{*}\sigma \left( P\right) =\sigma \left( P\right)\qquad
{\rm or}\qquad
 \sigma \left( P\right) +\alpha ^{*}\sigma \left( P\right) =1.
$$
\end{definition}

\del{For}\ins{To}
a spectral boundary value problem $\left(D,P\right) $ with a first-order
operator $D$ \ins{and} an even (odd) projection $P$\ins{,}
one can \del{define}\ins{assign}  a classical boundary
value problem. To this end, let us denote by $\alpha ^{*}D$ and
$\alpha ^{*}D^{-1}$ first-order elliptic operators with principal symbols
equal to $\alpha ^{*}\sigma \left(D\right) $ and
$\alpha ^{*}\sigma ^{-1}\left( D\right)$\ins{, respectively,}
on the sphere bundle $S^{*}M$.
In the even case\ins{,} the operator $D\oplus \alpha ^{*}D$ admits
\del{an}\ins{the}  elliptic classical boundary value problem
\begin{equation}
\left\{
\begin{array}{c}
Du=f_1,\qquad \alpha ^{*}Dv=f_2,
\vspace{0.2cm}\\
P\left. u\right| _X+\left( 1-P\right) \left. v\right| _X=g,\quad g\in
C^\infty \left( X,\left. E\right| _X\right) .
\end{array}
\right.   \label{a1}
\end{equation}
\del{Similarly}\ins{Likewise}, in the odd case
\del{there is a}\ins{we have the}  boundary value problem

\begin{equation}
\left\{
\begin{array}{c}
Du=f_1,\qquad \alpha ^{*}D^{-1}v=f_2,
\vspace{0.2cm}\\
P\left. u\right| _X+\left( 1-P\right) \left. v\right| _X=g,\quad g\in
C^\infty \left( X,\left. E\right| _X\right) .
\end{array}
\right.   \label{a2}
\end{equation}

\del{By}\ins{In the}  pass\del{ing}\ins{age}
f\del{o}r\ins{o}m the spectral
boundary value problem $\left( D,P\right) $
to the classical boundary value problem (\ref{a1}) or (\ref{a2}),
the dimension of the manifold $M$ \del{has to}\ins{must}
be taken into account.
The following proposition shows that \del{for}\ins{if}  the parity of
\ins{the} boundary value
problem\del{ which} is opposite to the parity of ${\rm
dim}M$\ins{,} \ins{then}
the boundary value problems (\ref{a1})\dell{,} \ins{and} (\ref{a2}) define
$2$-torsion elements in the group $\limfunc{Ell}^1\left( M\right)
\simeq K(T^*(M\setminus\partial M)) $.

\begin{proposition}
\label{utv}
\begin{enumerate}
\item  {\em The map\ins{ping} }$\alpha $
{\em induces \ins{an} involution in
}$K${\em -theory\ins{.} \del{that m}\ins{M}odulo
}$2${\em -torsion\ins{,} \ins{this involution is}
equal\del{s} \ins{to}  }$\left( -1\right) ^{\dim M}:$%
\[
\alpha ^{*}:K^{*}\left( T^{*}M\right) \otimes {\bf Z}\left[ \frac
12\right] \longrightarrow K^{*}\left( T^{*}M\right) \otimes
{\bf Z}\left[
\frac 12\right] ,\quad \alpha ^{*}=\left( -1\right) ^{\dim M}.
\]
{\em The involution }$\alpha ^{*}$ {\em has this property also
on the group }$K^{*}\left( T^{*}\left( M\backslash \partial
M\right) \right) .$

\item  {\em For an even-dimensional manifold }$M$, {\em
the projection }$S^{*}M\rightarrow P^{*}M ${\em
induces an isomorphism  (modulo }$2${\em-torsion) }
\[
K^{*}\left( P^{*}M\right) \otimes {\bf Z}\left[ \frac 12\right]
\rightarrow K^{*}\left( S^{*}M\right) \otimes {\bf Z}\left[ \frac
12\right]
\]
{\em (here }$P^{*}M=S^{*}M/\alpha$ {\em is the
corresponding projective cotangent sphere bundle).}

\item  {\em On an odd-dimensional manifold\ins{,} the projection }
$P^{*}M\rightarrow M$ {\em
induces an isomorphism (modulo }$2${\em-torsion) }
\[
K^{*}\left( M\right) \otimes {\bf Z}\left[ \frac 12\right] \rightarrow
K^{*}\left( P^{*}M\right) \otimes {\bf Z}\left[ \frac 12\right] .
\]
\end{enumerate}
\end{proposition}

\noindent {\em Proof.} Let us apply \ins{the} Mayer--Vietoris
principle \cite{BoTu1}.

\noindent 1) Let us check properties 1--3 for the restriction
of the mappings to the fiber over a point $x$ of the base
$M$ of the corresponding bundles:
\begin{eqnarray*}
K^{*}\left( T_x^{*}M\right) &\longrightarrow &K^{*}\left( T_x^{*}M\right) ,
\\
\,K^{*}\left( P_x^{*}M\right) &\rightarrow &K^{*}\left( S_x^{*}M\right) , \\
K^{*}\left( \left\{ x\right\} \right) &\rightarrow &K^{*}\left(
P_x^{*}M\right) .
\end{eqnarray*}
In the first case\ins{,} we have
\[
T_x^{*}M={\bf R}^{\dim M},\qquad K^{*}\left( {\bf R}^{\dim M}\right) =%
{\bf Z.}
\]
\ins{The} \del{I}\ins{i}nvolution $\alpha $
preserves (\ins{or} reverses) the orientation of the
space ${\bf R}^{\dim M}$\del{ together}
\del{with}\ins{depending on}  the parity of dimension of $M$.
Hence, we obtain the \del{requ}\ins{des}ired \ins{identity}
$\alpha^{*}=\left( -1\right) ^{\dim M}.$

In the second case\ins{,}
for an even-dimensional manifold $M$ we consider the
projection $\pi :S^{2n+1}\rightarrow {\bf RP}^{2n+1}$.
\ins{The} $K$-groups of
spheres and projective spaces are well-known (\del{see
}e.g.\ins{,} \ins{see} \cite{Gil11})\ins{:}
\begin{eqnarray*}
K^0\left( {\bf RP}^{2n+1}\right)  &=&{\bf Z\oplus Z}_{2^n},\quad
K^0\left( S^{2n+1}\right) ={\bf Z,} \\
K^1\left( {\bf RP}^{2n+1}\right)  &=&{\bf Z,\qquad }\quad K^1\left(
S^{2n+1}\right) ={\bf Z.}
\end{eqnarray*}
The first term in the groups $K^0$ is given by the dimension of vector
bundles, while the projection $\pi $ induces \del{a}\ins{the}
multiplication by 2
map\ins{ping}
on the groups $K^1$\ins{:}
\[
\begin{array}{ccc}
\pi ^{*}:K^1\left( {\bf RP}^{2n+1}\right) ={\bf Z} & {\bf
\longrightarrow } & K^1\left( S^{2n+1}\right) ={\bf Z},\\
%&  &  \\
n & \longmapsto  & 2n.
\end{array}
\]
In the third case\ins{,}
on an odd-dimensional $M$ we consider the projection
${\bf RP}^{2n}\rightarrow pt.$ The relevant $K$-groups are
\begin{eqnarray*}
K^0\left( {\bf RP}^{2n}\right)  &=&{\bf Z\oplus Z}_{2^n},\quad
K^0\left( pt\right) ={\bf Z,} \\
K^1\left( {\bf RP}^{2n+1}\right)  &=&0,\qquad K^1\left( pt\right) =0%
{\bf.}
\end{eqnarray*}
Both components ${\bf Z}$ correspond to
\ins{the} dimension of vector bundles.
Thus, \del{the }property $3$ is also satisfied over a point.

\noindent 2) By the Mayer--Vietoris principle\ins{,} we have to verify
the following assertion: \del{granted}\ins{if}  properties 1--3
are satisfied
over open subsets $U,V\subset M$ and\del{ on} their intersection $U\cap V,$
then these properties hold\del{ true} over the union $U\cup V$. %

$\!\!$In the first case\ins{,} let us write \del{down}\ins{out}
a part of \ins{the} Mayer--Vietoris
exact sequence
\[
\begin{array}{ccccc}
K^{*+1}\left( T^{*}\left( U\bigcap V\right) \right)  & \rightarrow  &
K^{*}\left( T^{*}\left( U\bigcup V\right) \right)  & \rightarrow  &
K^{*}\left( T^{*}U\right) \oplus K^{*}\left( T^{*}V\right)  \\
\downarrow \alpha ^{*} &  & \downarrow \alpha ^{*} &  & \quad \quad \quad
\downarrow \!\!\alpha ^{*}\!\!\oplus \!\alpha ^{*} \\
K^{*+1}\left( T^{*}\left( U\bigcap V\right) \right)  & \rightarrow  &
K^{*}\left( T^{*}\left( U\bigcup V\right) \right)  & \rightarrow  &
K^{*}\left( T^{*}U\right) \oplus K^{*}\left( T^{*}V\right) .
\end{array}
\]
\del{By a}\ins{A} diagram \del{search argument it can
be}\ins{chase}  show\del{n}\ins{s}  that the map\ins{ping}
$\alpha ^{*}$ in the middle satisfies property 1.

The second and the third cases \del{are}\ins{can be} treated
\ins{in a} similar\del{ly} \ins{way}.
For example, on an even-dimensional $M$\ins{,} the projection
$\pi :S^{*}M\rightarrow P^{*}M$ acts on the Mayer--Vietoris sequences
\[
\begin{array}{cccc}
\cdots \rightarrow \! & K^{*}\left( P^{*}\left( U\bigcup V\right) \right)%
\!\rightarrow \! & K^{*}\left( P^{*}\left( U\sqcup V\right) \right)
\rightarrow  & K^{*}\left( P^{*}\left( U\bigcap V\right) \right) \rightarrow
\cdots  \\ %
& \downarrow \pi ^{*} & \quad \downarrow \!\!\pi ^{*}\! & \downarrow \pi ^{*}
\\
\cdots \rightarrow \! & K^{*}\left( S^{*}\left( U\bigcup V\right) \right)%
\!\rightarrow \! & K^{*}\left( S^{*}\left( U\sqcup V\right) \right)
\rightarrow  & K^{*}\left( S^{*}\left( U\bigcap V\right) \right) \rightarrow
\cdots %
\end{array}
\]
By the \ins{five} \del{5-}lemma\ins{,}
the map\ins{ping} $\pi ^{*}$ on the left is an isomorphism
modulo 2-torsion.

The statement concerning the group
$K^{*}\left( T^{*}\left( M\backslash \partial M\right) \right) $ follows
from the exact sequence of the pair $\left. T^{*}M\right| _X\subset T^{*}M$
on which $\alpha^*$ acts\ins{:}
\[
\begin{array}{ccccccc}
\rightarrow\!\!& K^{*+1}\left( \left. T^{*}M\right| _X\right) \otimes
{\bf Z%
}\left[ \frac 12\right]  &\!\! \rightarrow\!\!  & K^{*}\left( T^{*}\left(
M\backslash \partial M\right) \right) \otimes {\bf Z}\left[ \frac
12\right]  & \rightarrow  & K^{*}\left( T^{*}M\right) \otimes {\bf Z}%
\left[ \frac 12\right]  &\!\! \rightarrow  \\
& \qquad \downarrow \left( -1\right) ^{\dim M} &  & \downarrow \alpha  &  &
\qquad \downarrow \left( -1\right) ^{\dim M} &  \\
\rightarrow\!\!& K^{*+1}\left( \left. T^{*}M\right| _X\right) \otimes
{\bf Z%
}\left[ \frac 12\right]  &\!\! \rightarrow\!\!  & K^{*}\left( T^{*}\left(
M\backslash \partial M\right) \right) \otimes {\bf Z}\left[ \frac
12\right]  & \rightarrow  & K^{*}\left( T^{*}M\right) \otimes {\bf Z}%
\left[ \frac 12\right]  &\!\! \rightarrow
\end{array}
\]

This completes the proof of Proposition \ref{utv}.

\vspace{0.1cm}

Thus, \del{below}\ins{in what follows}
we consider boundary value problems with even projections
$P$ on even-dimensional manifolds and odd projections on odd-dimensional
manifolds.

\subsection{The classification of
boundary value problems with even projections}

The boundary value problem\del{ in subspaces} (\ref{lastx})
\ins{in subspaces} can\del{
}not be
classified in terms of the classical boundary value problem
(\ref{a1}) or (\ref{a2}) even under the above\del{ introduced}
parity restrictions.
The point is that a classical boundary value problem is defined, up to
\ins{a} homotopy, by its principal symbol, while the boundary value problem
(\ref{lastx}) is not determined by the principal symbol. Indeed, by adding
\del{a }finite-dimensional space\ins{s}  to $\func{Im}P$, we obtain boundary
value problems with the same principal symbol\del{. The values
of the} \ins{but with different} index, \ins{which shows that}%
\del{ however, are different. Hence, the boundary value
problems} \ins{they} are
\ins{not} \del{non}homotopic \ins{to the original problem}.

It is shown in \cite{SaSt1,SaSt2} that the subspaces defined
by even (odd)
pseudodifferential projections have \del{a}\ins{the}
homotopy invariant\del{ that is}
described in the following theorem. Let us denote the
semigroups of subspaces
defined by even (odd) pseudodifferential projections by
$\widehat{\limfunc{Even}}\left(X\right) $ and $\widehat{\limfunc{Odd}}\left(
X\right)$\ins{,} \ins{respectively}.

\begin{theorem}
{\em \cite{SaSt1,SaSt2}}
There is a \ins{unique} homotopy invariant functional $d$
\[
d:\widehat{\limfunc{Even}}\left( X^{odd}\right) \rightarrow {\bf Z}\left[
\frac 12\right] ,\quad \text{or \quad }d:\widehat{\limfunc{Odd}}%
\left( X^{ev}\right) \rightarrow {\bf Z}\left[ \frac 12\right] ,
\]
\del{that is uniquely described by}\ins{with} the
\ins{following} properties:
\end{theorem}
\begin{enumerate}
\item  ({\em invariance\/})
\[
d\left( \func{Im}UPU^{-1}\right) =d\left( \func{Im}P\right)
\]
{\em for invertible pseudodifferential operators }$U${\em with
even principal symbol\/}: $\alpha ^{*}\sigma \left(
U\right) =\sigma \left( U\right)$;

\item  ({\em relative index\/})
\[
d\left( \func{Im}P_1\right) -d\left( \func{Im}P_2\right) =\limfunc{ind}%
\left( P_1,P_2\right),
\]
{\em where
${\rm ind}(P_1,P_2)={\rm ind}(P_2:\func{Im}P_1\to\func{Im}P_2)$
is the relative index of projections with equal principal symbols
}\cite{BDF1};

\item  ({\em complement\/})
\[
d\left( \func{Im}P\right) +d\left( \func{Im}\left( 1-P\right) \right) =0.
\]
\end{enumerate}

The group of stably homotopic spectral boundary value problems with even
projections $P$ is denoted by $\limfunc{Ell}^{ev}\left( M,
\partial M\right).$ It \del{appears}\ins{turns out}
that the classical boundary value problem
(\ref{a1})\del{ together} \del{with}\ins{and}  the invariant $d$ of the subspace of right-hand
sides \ins{already}
classify spectral boundary value problems modulo $2$-torsion.

\begin{theorem}
\label{theven}On an even-dimensional manifold $M$\ins{,}
the map\ins{ping}
\[
\begin{array}{ccc}
\chi :\limfunc{Ell}^{ev}\left( M,\partial M\right) \otimes {\bf Z%%
}\left[ \frac 12\right]  & \rightarrow  & \limfunc{Ell}^1\left(
M\right) \otimes {\bf Z}\left[ \frac 12\right] \oplus {\bf Z}\left[
\frac 12\right] , \\
&  &  \\
\left( D,P\right)  & \longmapsto  & \left( \left( D\oplus \alpha
^{*}D\right) \otimes \frac 12,d\left( \func{Im}P\right) \right)
\end{array}
\]
is an isomorphism of abelian groups.
\end{theorem}

\noindent {\em Proof. }Let us define the inverse map\ins{ping}
\[
\chi ^{\prime }:\limfunc{Ell}^1\left( M\right) \otimes {\bf Z}%%
\left[ \frac 12\right] \oplus {\bf Z}\left[ \frac 12\right] \rightarrow
\limfunc{Ell}^{ev}\left( M,\partial M\right) \otimes {\bf Z}%
\left[ \frac 12\right] .
\]
On the first term it is induced by the embedding of classical boundary value
problems in\del{to}
boundary value problems with even projections, while on the
second term it is given by the formula
\[
\chi ^{\prime }\left( 0,\frac k{2^N}\right) =\frac 1{2^N}\left(%
D_{+},k\right) ,
\]
where $\left( D_{+},k\right) $
\del{denotes}\ins{stands for}  the spectral boundary value problem
for the operator $D_{+}$ with a finite-dimensional spectral projection of
rank $k.$ Let us verify that $\chi ^{\prime}$ is the inverse
\del{for}\ins{of}  $\chi.$

The second component of the composition $\chi \circ\chi ^{\prime }$
\[
{\bf Z}\left[ \frac 12\right] \rightarrow {\bf Z}\left[ \frac
12\right]
\]
is the identity map\ins{ping}
\del{due to}\ins{by}  property 2) of the functional $d.$ The first
component \ins{is} equal\del{s} \ins{to}
\[
\frac{\alpha ^{*}+1}2:\limfunc{Ell}^1\left( M\right) \otimes
{\bf Z}\left[ \frac 12\right] \longrightarrow \limfunc{Ell}%
^1\left( M\right) \otimes {\bf Z}\left[ \frac 12\right] ,
\]
\del{that}\ins{which}, by virtue of \del{an}\ins{the}  isomorphism
\[
\limfunc{Ell}^1\left( M\right) \simeq K\left( T^{*}\left(
M\backslash \partial M\right) \right)
\]
and Proposition \ref{utv}, \del{p.}\ins{item}~1,
is the identity map\ins{ping}.

The assertion of the theorem \ins{can} now
\del{follows}\ins{be derived}  from the \ins{following} lemma.

\begin{lemma}
\label{hitro}The homomorphism $\chi ^{\prime }$ is an epimorphism.
\end{lemma}

\noindent {\em Proof. }Consider an arbitrary spectral boundary value
problem $\left(D,P\right) $\dell{,} with
even projection $P$. Proposition \ref{utv},
\del{p.}\ins{item~}3 implies that the \ins{sum of $2^N$ copies
of the}
subbundle $\func{Im}\sigma \left( P\right) \subset
\pi ^{*}E$\del{ being taken $2^N$ copies}
is homotopic in the class of even
subbundles
to a bundle $\pi ^{*}G,$ $G\subset E^{2^N}$,
\del{pulled-back}\ins{lifted} from the base $X$.
We denote the corresponding homotopy of projections by
$\sigma \left(P_t\right)$:
\[
\sigma \left( P_t\right) :\pi ^{*}E^{2^N}\longrightarrow \pi
^{*}E^{2^N},\qquad t\in \left[ 0,1\right] ,
\]
\[
\sigma \left( P_0\right) =2^N\sigma \left( P\right) ,\quad \func{Im}\sigma
\left( P_1\right) =\pi ^{*}G.
\]
Consider \del{the}\ins{a}
covering homotopy of pseudodifferential projections
$P_t$\dell{,}
such that $P_0=2^NP.$ The symbol of\del{ projection} $P_1$ is
equal to the symbol of projection $P_G$ on\del{to} the space
$C^\infty \left( X,G\right) \subset C^\infty \left( X,E\right) $
of sections of bundle $G.$ Hence, the homotopy classification
of projections with the same principal symbols \cite{Wojc1} shows
that $P_1$ is homotopic to a projection differing from $P_G$
by a finite\dell{-}\ins{ }rank projection. We can assume that the homotopy
$P_t$ already gives\del{ for $t=1$} such \ins{a} projection
\ins{at $t=1$}.

The homotopy of projections $P_t$ extends to a homotopy of spectral
boundary value problems
\[
\left( D_t,P_t\right) :D_0=2^ND
\]
by\del{ the} formula (\ref{dva}).
The spectral boundary value problem
$\left(D_1,P_1\right) $ then lies in the range of the
map\ins{ping} $\chi^{\prime }$ \ins{given by}
\[
\chi ^{\prime }\left[ \left( D_1,P_G\right) ,-\limfunc{ind}\left(%
P_G,P_1\right) \right] =\left[ D_1,P_1\right] .
\]
This proves the lemma.
The theorem is \ins{thereby} proved.

\vspace{0.1cm}

\subsection{The classification of
boundary value problems with odd
projections}

Let us generalize the definition of spectral boundary value problems with
odd projections. We consider spectral boundary value
problems\dell{,} such that
the symbol of \ins{the} projection $P$ is \del{a}\ins{the}
sum of a constant symbol with respect
to the cotangent variables and an odd projection.  Let us also
identify spectral boundary value problems of the form
\begin{equation}
\label{lstar}
\left\{
\begin{array}{c}
D_1u=f_1,\qquad D_2v=f_2, \\
P\left. u\right| _X=g_1,\quad \left( 1-P\right) \left. v\right| _X=g_{2}
\end{array}
\right.
\end{equation}
with odd projection $P$ with the corresponding classical boundary
value problems
\begin{equation}
\label{lstar1}
\left\{
\begin{array}{c}
D_1u=f_1,\qquad D_2v=f_2, \\
P\left. u\right| _X+\left( 1-P\right) \left. v\right| _X=g.
\end{array}
\right.
\end{equation}
\del{An}\ins{The}
abelian group of stable homotopy classes of such spectral boundary value
problems \del{is}\ins{will be}
denoted by $\limfunc{Ell}^{odd}\left( M,\partial M\right)$.
In the following theorem\ins{,} the stable homotopy classification modulo
2-torsion is established for spectral boundary value problems with
odd projections.

\begin{theorem}
\label{thodd}On an odd-dimensional manifold $M$\ins{,} the
map\ins{ping}
\[
\begin{array}{ccc}
\chi :\limfunc{Ell}^{odd}\left( M,\partial M\right) \otimes {\bf
Z}\left[ \frac 12\right]  & \rightarrow  & \limfunc{Ell}^1\left(
M\right) \otimes {\bf Z}\left[ \frac 12\right] \oplus {\bf Z}\left[
\frac 12\right] , \\
&  &  \\
\left( D,P\right)  & \longmapsto  & \left( \left( D\oplus \alpha
^{*}D^{-1}\right) \otimes \frac 12,d\left( \func{Im}P\right) \right)
\end{array}
\]
is an isomorphism of abelian groups.
\end{theorem}

\noindent {\em Proof. }Let us define the inverse map\ins{ping}
\[
\chi ^{\prime }:\limfunc{Ell}^1\left( M\right) \otimes {\bf Z}%%
\left[ \frac 12\right] \oplus {\bf Z}\left[ \frac 12\right] \rightarrow
\limfunc{Ell}^{odd}\left( M,\partial M\right) \otimes {\bf Z}%
\left[ \frac 12\right] .
\]
On the first summand it is induced by the embedding of classical
boundary value problems in the class of boundary value problems with
odd projections\ins{,} and on the second summand it is given,
\ins{just} as in the previous
theorem, by the formula
\[
\chi ^{\prime }\left( 0,\frac k{2^N}\right) =\frac 1{2^N}\left(%
D_{+},k\right) .
\]

The second component of the composition $\chi \circ\chi ^{\prime }$
\[
{\bf Z}\left[ \frac 12\right] \rightarrow {\bf Z}\left[ \frac
12\right]
\]
is equal to the identity map\ins{ping}, while the first component is
\[
\frac{1-\alpha ^{*}}2:\limfunc{Ell}^1\left( M\right) \otimes
{\bf Z}\left[ \frac 12\right] \longrightarrow \limfunc{Ell}%
^1\left( M\right) \otimes {\bf Z}\left[ \frac 12\right],
\]
and, by virtue of the isomorphism
\[
\limfunc{Ell}^1\left( M\right) \simeq K\left( T^{*}\left(
M\backslash \partial M\right) \right)
\]
and Proposition \ref{utv}, \del{p.}\ins{item} 1, it is equal to
\ins{the} identity \ins{mapping}.

The following lemma completes the proof of the theorem.

\begin{lemma}
\ins{The} \del{H}\ins{h}omomorphism $\chi ^{\prime }$ is an epimorphism.
\end{lemma}

\noindent {\em Proof. }Consider the spectral boundary value problem
$\left(D,P\right) $ with an odd projection $P.$
\ins{Just} \del{A}\ins{a}s in the proof of\del{ previous} Lemma \ref{hitro},
we only need to construct a homotopy of the principal symbol
of\del{ projection}
$P$ to a projection independent of \ins{the} cotangent variables.
By virtue of the identification (\ref{lstar}),
(\ref{lstar1})\ins{,}
it\del{ is} suffic\del{ient}\ins{es}
to construct a homotopy of the principal symbol
of \ins{the} projection to a direct sum
$p\oplus \left(1-p\right)$, where $p$ is an odd projection.

It \del{is}\ins{was}
proved in \cite{SaSt2} that for some $N$ there exists an even
isomorphism
\[
u:2^N\pi ^{*}\left. E\right| _X\rightarrow 2^N\pi ^{*}\left. E\right|
_X,\qquad \pi :S^{*}X\rightarrow X,
\]
that \del{transforms}\ins{takes}
the projection $2^N\sigma \left( P\right) $ \del{in}to the
complementary projection $2^N\left( 1-\sigma \left( P\right)
\right)$.
Moreover, this isomorphism defines \ins{the} zero element in the group
$K^1\left( S^{*}X\right)$.
\ins{It follows} \del{F}\ins{f}rom Proposition \ref{utv}\ins{,}
\del{p.}\ins{item~}2\del{ it follows} that for
\ins{sufficiently} large\del{ enough}
$N$ the isomorphism $u$ is homotopic to \ins{the}
identity in the class of even
isomorphisms. Let us denote a homotopy of this type by $u_t$,
$t\in\left[ 0,1\right]$, $u_0=1.$ The
\del{required}\ins{desired}  homotopy of projections
is \del{defined}\ins{given} by the formula
\[
\sigma \left( P_t\right) =2^N\sigma \left( P\right) \oplus u_t2^N\sigma
\left( P\right) u_t^{-1},
\]
\[
\sigma \left( P_0\right) =2^{N+1}\sigma \left( P\right) ,\qquad \sigma
\left( P_1\right) =2^N\left[ \sigma \left( P\right) \oplus \left( 1-\sigma
\left( P\right) \right) \right].
\]

\ins{The} \del{L}\ins{l}emma\del{ is proved.} \ins{and}
\del{The}\ins{the}  theorem \del{is}\ins{are thereby}  proved.

\begin{corollary}
{\em\cite{SaSt1,SaSt2}} Spectral boundary
value problems $\limfunc{Ell}^{ev/odd}\left(M,\partial M\right) $
with parity con\-ditions have the following homotopy classification
modulo 2-torsion\ins{:}
$$
\limfunc{Ell}^{ev/odd}\left(M,\partial M\right)
     \otimes{\bf Z}
     \left[ \frac 12\right]
  \simeq
K\left(T^*(M\setminus\partial M)\right) \otimes {\bf Z}%%
\left[ \frac 12\right] \oplus {\bf Z}\left[ \frac 12\right].
$$
The index formula\del{ is valid}
\begin{equation}
\limfunc{ind}\left( D,P\right) =\frac 12\limfunc{ind}\left( D\oplus \alpha
^{*}D^{\pm 1}\right) -d\left( \func{Im}P\right)   \label{quest}
\end{equation}
\ins{is valid.}
\end{corollary}

\noindent Indeed, let us consider both sides of the index formula as
homomorphisms of the group $\limfunc{Ell}^{ev/odd}\left(M,\partial M\right) $
in\ins{to}  ${\bf Z}\left[ \frac 12\right]$.
By Theorems~\ref{theven} and~\ref{thodd}, the groups
$\limfunc{Ell}^{ev/odd}\left(M,\partial M\right) $ are
rationally generated by classical boundary value problems and\del{ also
by}
boundary value problems with finite-dimensional spectral
subspaces. On both types of generators\ins{,}
the two parts of the index formula
(\ref{quest}) coincide. This proves the index formula.

%\bibliography{elliptic}
%\bibliographystyle{unsrt1}
\vspace{1cm}

\hfill {\em Moscow, Potsdam}
\end{document}